%-----------------------------------------------
%         GABRIELOV'S EXAMPLE   HH   JULY 2015
%           APRIL/AUG 2016 JUNE 2017 FEB 2018
%-----------------------------------------------

\input amssym.tex
\parindent=0cm
\hoffset=0.8truecm
\hsize=14.5truecm
\vsize=25truecm
%\hsize=16.5truecm
%\hoffset=-0.2truecm
%\vsize =25.5cm
\voffset=-.5cm
\baselineskip 15pt
%\nopagenumbers
\long\def\ignore#1\recognize{}

%-----------------------------------------------
%    FONTNAMES  FOR WINEDT  + TEXSHOP
%-----------------------------------------------

%\ignore

\font\Times=ptmr at 10pt
\font\Bf=ptmb at 12pt
\font\bf=ptmb at 10pt

\font\smallTimes=ptmr at 8pt
\font\smallit=ptmri at 8pt

\Times

%\recognize
%-----------------------------------------------
%                MACROS
%-----------------------------------------------

\def\litem{\par\noindent\hangindent=\parindent\ltextindent}
\def\ltextindent#1{\hbox to \hangindent{#1\hss}\ignorespaces}
\long\def\ignore#1\recognize{}

\def\big{\bigskip}
\def\med{\medskip}
\def\gb{\goodbreak}
\def\qed{\hfill$\circlearrowleft$\med}
\def\hs{\hskip}
\def\vs{\vskip}

\def\cl{\centerline}
\def\ds{\displaystyle}
\def\ol{\overline}

\def\wh{\widehat}
\def\wt{\widetilde}
\def\sm{\setminus}

\def\{{\lbrace}
\def\}{\rbrace}

\def\map{\rightarrow}
\def\inv{^{-1}}
\def\N{{\Bbb N}}

\def\C{{\Bbb C}}
\def\Z{{\Bbb Z}}

\def\Q{{\Bbb Q}}

\def\abs#1{\vert#1\vert}
\def\<{<\!\!}
\def\>{\!\!>}

\def\a{\alpha}
\def\b{\beta}
\def\6{\partial}
\def\min{{\rm min}}
\def\max{{\rm max}}

\def\Id{{\rm Id}}
\def\inin{{\rm in}}
\def\int{{\rm int}}

\def\supp{{\rm supp}}

\def\.{$\bullet$}

\def\eps{{}}%{\varepsilon}
%-------------------------------------------------

\def\hh{$\bullet$ }

%\ignore

%-----------------------------------------------
%         GABRIELOV'S EXAMPLE   HH   JULY 2015 APRIL/AUG 2016 JULY 2017
%-----------------------------------------------

\cl{\Bf Echelons of power series and Gabrielov's counterexample}\med

\cl{\Bf  to nested linear Artin Approximation}
\big\big

%-------------------------------------------------

\cl{{M.E. ALONSO, F.J. CASTRO-JIM\'ENEZ, H. HAUSER, C. KOUTSCHAN}
%
%-------------------------------------------------
\footnote{}{\smallTimes Mathematics Subject Classification (2010): 32A05, 13F25, 13P10, 16W60, 14QXX, 14B12. }
%-------------------------------------------------
\footnote{${}^{(1)}$}{\smallTimes This work was done in part during a Research-in-Teams program at the Erwin-Schr\"odinger Institute at Vienna, and the special semester on Artin approximation within the Chaire Jean Morlet at CIRM, Luminy-Marseille. M.E.A. was supported by MINECO MTM2014-55565  and   UCM , IMI and Grupo 910444, F.J.C.-J. by MTM2013-40455-P and MTM2016- 75024-P, H.H. and C.K. by the Austrian Science Fund FWF, within the projects P-25652 and AI-0038211, respectively P29467-N32 and F5011-N15.
}}
%-------------------------------------------------
\big\big

\cl{\vbox {\hsize 10 cm \smallTimes {\smallit Abstract} : Gabrielov's famous example for the failure of analytic Artin approximation in the presence of nested subring conditions is shown to be due to a growth phenomenon in standard basis computations for {\smallit echelons}, a generalization of the concept of ideals in power series rings.}}\par

%\vs.5cm\hfill hh aug 4, 2017
\big

%-------------------------------------------------

{\bf Introduction}\med

In the S\'eminaire Henri Cartan of 1960/61, Grothendieck posed the question whether analytically independent analytic functions are also formally independent [Gr].\footnote
%
%-------------------------------------------------
{${}^{(2)}$}{\smallTimes Artin attributes in [Ar] Grothendieck's question to Abhyankar.}
%-------------------------------------------------
%
%[More specifically, one can ask whether the ideal $\wh R$ of formal relations $\wh r(y)\in\C[[y_1,\ldots,y_m]]$ between given analytic functions $f_1(x),\ldots,f_m(x)$ in variables $x_1,\ldots,x_n$, say,\med
%
%\hs 3cm  $\wh r(f_1(x),\ldots,f_m(x)) =0$, \med
%
%is already generated by the analytic relations $r(y)\in\C\{y_1,\ldots,y_m\}$.]
%
It came as a surprise when Gabrielov answered the question in 1971 in the negative. He constructed four analytic functions $e,f, g,h$ in three variables admitting one formal relation but no analytic one [Gb1]. To our knowledge, this is essentially the only known counterexample to Grothendieck's question. In an opposite direction, Paw\l ucki constructed analytic functions and a subset $Z$ of the reals for which there do exist analytic relations for parameter values outside $Z$ but there do not exist formal relations for parameters in $Z$ [Pa1]. In a later paper, Gabrielov gave a sufficient condition for a positive answer to Grothendieck's question in terms of the rank of the Jacobian matrix of the analytic functions [Gb2], see also [Pa2]. %The associated ring maps are nowadays called "regular" morphisms. %[see in this context Cutkosky's paper [Cu] achieving regularity via blowups, and Pfister's write-up from Luminy].
Much more generally, Popescu proved in 1985 a difficult approximation theorem which contains as a particular case a positive answer whenever the analytic functions are {\it algebraic} power series [Po1, Po2, Sp, Te]. %\med
Gabrielov's counterexample is based on an example of Osgood [Os] from 1916, complemented by a tricky construction and calculation. The deeper reason for the existence of formal divergent relations between analytically independent analytic functions remained mysterious over the years.\med

In this note we explain the genesis of the phenomenon in Gabrielov's example and provide a systematic way to construct many more counter\-examples: It turns out that the existence of formal but not analytic relations is caused by {\it accumulated growth occurrences} in standard basis computations for echelons (an {\it echelon} is a generalization of an ideal in a power series ring, see below). Such a growth behaviour is well known for standard bases of ideals, but does not do any harm there due to the finiteness of the basis (which is ensured by the Noetherianity of the power series ring.) Standard bases of echelons need no longer be finite, and the iterated growth occurrence in their construction may indeed force divergence. We illustrate in the paper how this phenomenon is related to the presence of sufficiently fast converging coefficients of the (analytic) input series. In the example, the coefficients converge faster than exponentially.\med

For algebraic power series, the phenomenon does not happen. The echelon standard basis may still be infinite, but the convergence of the coefficients of the involved series seems to be sufficiently slow so as to ensure a positive answer to Grothendieck's question: whenever there is a formal linear relation respecting the scopes, there is also a convergent one (actually, even an algebraic one). The assertion for algebraic series follows for instance from Popescu's approximation theorem (i.e., the fact that nested approximation holds for algebraic power series), whereas a direct explanation in terms of echelons is still lacking. \med

Our explanation of Gabrielov's example will be embedded in a short description of the division theorem for power series in the setting of echelons and the related notion of echelon standard basis. This is not mandatory to understand the example but should allow the reader to see its construction in a broader context.
\med\goodbreak
%-------------------------------------------------
%          GABRIELOV'S EXAMPLE INTRO
%-------------------------------------------------

{\bf Gabrielov's example}\med

The first step towards Grothendieck's question, and this already appears in [Gb1], is to transcribe the existence of formal or analytic relations to a nested linear Artin
approximation problem: Let $f_1(x),\ldots,f_m(x)$ be convergent power series in variables $x=(x_1,\ldots,x_n)$ and let $r(y_1,\ldots,y_m)$ be a (formal or analytic) relation between them, say,\med

\hs 3cm     $r(f_1(x),\ldots,f_m(x)) = 0$.\med

This is equivalent to saying that $r(y)$ belongs to the ideal of the formal, respectively convergent, power series ring $\C[[x,y]]$, respectively $\C\{x,y\}$, generated by the series $y_i - f_i(x)$, for $i=1,\ldots,m$. Therefore there exist power series $a_1(x,y),\ldots,a_m(x,y)$ such that\med

\hs 3cm   $\ds r(y) = \sum_{i=1}^m a_i(x,y) \cdot (y_i - f_i(x))$.\med

Here, the series $a_i$ are allowed to depend on both $x$ and $y$, whereas the series $r$ must be independent of $x$. This requirement is known in the context of  Artin approximation as a ``nested subring condition''. Note that the unknown series $r$ and $a_i$ appear linearly in the equation. As an extension of Grothendieck's question one may then ask more generally whether linear nested Artin approximation holds for analytic functions: Given analytic functions $e$ and $f_1,\ldots,f_m$  in $n$ variables $x_1,\ldots,x_n$ such that the linear presentation\med

\hs 3cm $\ds e(x)=\sum_{i=1}^m \wh a_i(x)\cdot f_i(x)$\med

holds with formal power series $\wh a_i(x)$ depending only on the variables $x_1,\ldots,x_{s_i}$, for given $s_i\leq n$, does there exist a presentation \med

\hs 3cm $\ds e(x)=\sum_{i=1}^m a_i(x)\cdot f_i(x)$\med

with analytic functions $a_i(x)$ depending on the same sets of variables as $\wh a_i(x)$?\med\goodbreak

% A difficult theorem of Popescu asserts that if the series $f_i$ are algebraic power series (i.e., formal power series which are algebraic over the field $\C(x)$ of rational functions), then the existence of a formal nested solution $\wh r(y), \wh a_1(x,y),\ldots,\wh a_m(x,y)$ to the equation $(*)$ already implies the existence of an algebraic (and hence, in particular, convergent) solution $r(y), a_1(x,y),\ldots,a_m(x,y)$ [Po1, Po2, Sp, KPPRM]. Gabrielov's example shows that this is false if the series $f_i(x)$ are not algebraic but merely convergent.
\med
\goodbreak
%-------------------------------------------------

Gabrielov also gives a counterexample to this case of linear nested analytic approximation: Consider the series $f = 1$, $g = x\cdot (e^z-1)$,  and $h = yz - x$ in three variables $x, y, z$. He then shows that the convergent series
%
%-------------------------------------------------
%\footnote{${}^{(3)}$}{Actually, Gabrielov takes for $g$ the series $g'=zy\cdot e^z$. Then the linear combination $g'- z^e\cdot h - x\cdot f=zy\cdot e^z- e^z\cdot (zy-x) - x= x\cdot (e^z-1)$ gives our $g$. He then works with the series $e'(y,z)=\sum_{i=1}^\infty\sum_{j=0}^\infty {i!\over (i+j)!}\cdot y^iz^{i+j+1} $ which, modulo $h$, is equivalent to the series $e(x,z)$ used above.}
%-------------------------------------------------
%
 \med

\hs3cm  $\ds e(x,z)=\sum_{i=1}^\infty\sum_{j=0}^\infty {i!\over (i+j)!}\cdot x^iz^{j+1}$\med

admits a presentation \med

\hs 3cm $e=\wh a\cdot f+\wh b\cdot g+\wh c\cdot h$,\med

with formal series $\wh a(x,y)$, $\wh b(x,y)$, $\wh c(x,y,z)$ but that there are no convergent series $a(x,y)$, $b(x,y)$, $c(x,y,z)$ representing $e$ in this way. Setting\med

\hs3cm $ I = \C\{x,y\}\cdot f+ \C\{x,y\}\cdot g + \C\{x,y,z\}\cdot h$\med

with completion $\wh I=\C[[x,y]]\cdot f+ \C[[x,y]]\cdot g + \C[[x,y,z]]\cdot h$, one therefore has the strict inclusion of vector subspaces\med

\hs 3cm $I\subsetneq \wh I\cap \C\{x,y,z\}$.\med

We will investigate in this note the deeper reason behind this fact. To do so, we collect  in the next section the basics about echelons, a generalization of the notion of ideals in power series rings. Subspaces as $I$ above are echelons, and their understanding is crucial for explaining Gabrielov's example. This explanation is presented in the section after the section on echelons.\big

\goodbreak
%-------------------------------------------------
%            ECHELONS
%-------------------------------------------------

{\bf Echelons}\med

%-------------------------------------------------
%              DEF ECHELONS + SCOPES
%-------------------------------------------------

Let $x_1,\ldots,x_n$ be variables and $K[[x]]=K[[x_1,\ldots,x_n]]$ the ring of formal power series in $x_1,\ldots,x_n$ over a given field $K$. If $K$ is a valued field, we may also consider the subring $K\{x\}=K\{x_1,\ldots,x_n\}$ of convergent power series. The next definitions apply always equally to the convergent case. A {\it finitary echelon} is a $K$-subspace of $K[[x]]$ which can be written as a finite sum \med

\hs 3cm $\ds I=\sum_{i=1}^k\, K[[x_1,\ldots,x_{s_i}]]\cdot f_i$,\med

with series $f_i\in K[[x_1,\ldots,x_n]]$ and integers $0\leq s_i\leq n$, called the {\it assigned scope} of $f_i$. Sometimes, we also refer to the variables $x_1,...,x_{s_i}$ themselves as the scope of $f_i$. Series $f_1,...,f_k$ as above with assigned scopes $s_1,...,s_k$ are called {\it generators} of $I$.  When working with finitary echelons, we often tacitly assume that a generator system is already chosen. A linear combination $f=\sum_{i=1}^k a_i\cdot f_i$ is said {\it to respect the scopes} if $a_i\in K[[x_1,...,x_{s_i}]]$ holds for all $i$. Each element $f$ of $I$ can be represented in this way. In certain situations, the sum $\sum_{i=1}^k\, K[[x_1,\ldots,x_{s_i}]]\cdot f_i$ will be direct, and then the presentation $f=\sum_{i=1}^k a_i\cdot f_i$ of elements $f\in I$ as a linear combination of $f_1,...,f_k$ respecting the scopes is unique (compare this with the later analysis of Gabrielov's example where the involved echelon is indeed a direct sum). \med

One could also develop a concept of infinitely generated echelons, but this is not needed for the sequel, and will hence be omitted here.\med

For $f\in I$, we call $s(f)=\max\{s\in\{0,\ldots,n\},\, K[[x_1,\ldots,x_s]]\cdot f\subset I\}$ the {\it actual scope} of $f$ in $I$. Clearly, the assigned scope $s$ of $f$ is less than or equal to the actual scope $s(f)$. In a theoretical context, we may always assign the actual scope to $f$, so that $s=s(f)$, and we then just speak of the scope of an element. But for actual computations and a given $f\in I$, it seems often impossible to determine the actual scope algorithmically, since it would require a constructive echelon membership test for the multiples of $f$. This aspect will play a role in Thm.~2, where only assigned scopes are considered.\med

%-------------------------------------------------
%         MODULE ECHELONS
%-------------------------------------------------

The analogous definitions hold for $K$-subspaces of free modules $K[[x]]^m$ of the form\med

\hs 3cm $\ds I=\sum_{i=1}^k\, K[[x_1,\ldots,x_{s_i}]]\cdot f_i\subset K[[x]]^m$,\med

with power series vectors $f_i\in K[[x_1,\ldots,x_n]]^m$. We call such subspaces {\it finitary (module) echelons}, with generators $f_i$ and assigned scopes $s_i$.\med

Let now $f_1,...,f_k$ with scopes $s_1,...,s_k$ be given generators of a finitary echelon $I\subset K[[x]]$ (or of a finitary module echelon $I\subset K[[x]]^m$). The {\it (module) echelon of (linear) relations} between $f_1,...,f_k$ is the $K$-subspace\med

\cl{${\rm Rel}(f_1,...,f_k)= \{r\in \prod_{i=1}^k K[[x_1,...,x_{s_i}]],\, \sum r_if_i=0\}$}\med

of $K[[x]]^k$ consisting of the linear relations between $f_1,...,f_k$ respecting the scopes $s_i$. Here, we assign to a relation $r=(r_1,...,r_k)$ the scope $t:=\min\{s_i,\, r_i\neq 0\}$, so that the inclusion $K[[x_1,...,x_t]]\cdot r\subset {\rm Rel}(f_1,...,f_k)$ is ensured. \med

%-------------------------------------------------
%      RELATIONS FINITELY GENERATED
%-------------------------------------------------

\ignore

We claim that ${\rm Rel}(f_1,...,f_k)$ is again a finitary (module) echelon, i.e., that there exist relations $r_j$ with scope $t_j$, for $j=1,...,p$, so that\med

\cl{${\rm Rel}(f_1,...,f_k)= \sum_{j=1}^p\, K[[x_1,\ldots,x_{t_j}]]\cdot r_j$.}\med

[\hh Please comment what to do with this: We should include a proof of this, though we will not use the statement except in the trivial case where the $f_i$ are monomials. But then one could talk about free resolutions of echelons and try to prove that every finitary echelon admits a finite free resolution. This could be an interesting result. For the proof of the finite generation of the relation echelon, one could probably use division, but then we would need to formulate the division theorem for module echelons: It is clear that the relation echelon between monomials $f_i$ (or monomial vectors having only one non-zero entry) is finitary, and then, for arbitrary $f_i$, use that, by the division theorem, an inclusion of echelons is an equality if their initial echelons coincide.]\med

\recognize
%-------------------------------------------------
%           MONOMIAL ORDERS
%-------------------------------------------------

Assume that a monomial order $<_\eps$ on $\N^n$ is chosen, i.e., a total ordering compatible with the addition in $\N^n$ and so that $0$ is the smallest element. It induces an ordering, also denoted by $<_\eps$, on the set of monomials $x^\a$ of $K[[x]]$, $\a\in\N^n$. For $f\in K[[x]]$, we denote by $\inin_\eps(f)=\inin(f)=x^\a$ the smallest monomial with respect to $<_\eps$ appearing in the expansion of $f$ (we always take the coefficient equal to $1$, and agree that $\inin(0)=0$). It is called the {\it initial monomial} of $f$ with respect to $<_\eps$. For a finitary echelon $I$ in $K[[x]]$, denote by $\inin_\eps(I)=\inin(I)$ the associated {\it initial echelon} of $I$: this is the $K$-subspace of $K[[x]]$ of power series whose expansion only involves monomials which are initial monomials $\inin(f)$ of elements $f$ of $I$. It is thus the $x$-adic closure of the subspace of $K[[x]]$ spanned by all initial monomials of elements of $I$. In general, $\inin(I)$ will not be a {\it finitary} echelon. %\hh [Similar definitions can be made for power series vectors in $K[[x]]^m$ with respect to (suitably defined) monomial orders $<_\eps$ on $\N^n\times\{1,...,m\}$, yielding initial monomial vectors $\inin_\eps(f)\in K[x]^m$ which have just one non-zero entry, and this entry is a monomial.]
For later use we restrict to monomial orders which admit no infinite bounded and strictly increasing sequences (thus, $(\N^n,<)$ will be order equivalent to $\N$ with the usual order). We reserve the symbol (\#) for this condition; it would not hold for instance for a lexicographic monomial order on $\N^n$. \med

Assume that we are given generators $F=\{f_1,\ldots,f_k\}$ of $I$ with assigned scopes $s_1,\ldots,s_k\in\{0,\ldots,n\}$,\med

\cl{$ I =\sum_{i=1}^k\, K[[x_1,\ldots,x_{s_i}]]\cdot f_i$.}\med

%The {\it actual scope} of an element $f$ of $I$ is the maximal number $s=s(f)\leq n$ so that $K[[x_1,\ldots,x_s]]\cdot f\subset I$. Clearly, one has $s_i\leq s(f_i)$. Here, we sometimes just say scope instead of actual scope, whereas assigned scopes are always called like this.

%Note that for the generators $f_1,\ldots,f_k$ of $I$ we may assume from the beginning that the assigned scope is already the actual scope. %Later on it might be necessary to distinguish between the two notions [namely, assigned scopes will be used to define decompositions of finitary echelons into {\it direct} sums of slices or azulejos].\med

Our goal is to construct an {\it echelon standard basis} of $I$ from $f_1,\ldots,f_k$. This is a (possibly infinite) set of elements $g_j$ of $I$, $j\in \N$, with assigned scopes $t_j$, whose initial monomials $x^{\a_j}=\inin(g_j)$ generate $\inin(I)$ topologically:\med

\cl{$ \inin(I) = \sum^*_{j\in\N}\, K[[x_1,\ldots,x_{t_j}]]\cdot x^{\a_j}$.}\med

Here, the symbol $\sum^*$ stands for {\it infinite} sums of elements of the summands $K[[x_1,\ldots,x_{t_j}]]\cdot x^{\a_j}$, say, power series in $K[[x]]$ whose exponents belong to the set \med

\cl{$\bigcup_{j\in\N}\, \a_j+(\N^{t_j}\times 0^{n-t_j})$.}\med

Such sums converge in the $x$-adic topology of $K[[x]]$ since the total degree of the monomials $x^{\a_j}$ tends with $j$ towards infinity (here, we exclude wlog repetitions among these monomials). By ``construction'' we understand a possibly infinite algorithm, which ``terminates'' in the sense that it produces, for each initial monomial $x^\a$ of $\inin(I)$, in finitely many steps an element $f\in I$ together with an assigned scope $s$ such that $x^\a\in K[[x_1,\ldots,x_{s}]]\cdot \inin(f)$. Our algorithm mimicks Buchberger's algorithm for the construction of Gr\"obner and/or standard bases of ideals of polynomials, respectively power series [Bu, GP]. We do not, however, divide the new elements after each step by the existing ones. The main difference to the case of ideals is that echelon standard bases are not necessarily finite sets of generators, so that the notion of ``termination'' of the algorithm has to be drafted properly. See Thm.~2 below for details.\footnote
{${}^{(3)}$}{\smallTimes In the case of {\smallit ideals} of power series rings, standard bases are finite. If the ideals are generated by algebraic power series (and the coordinates are sufficiently generic), the construction of standard bases can be performed by a finite algorithm, see [ACH].}\med

%-------------------------------------------------

Denote by $x^{\a_i}$ the initial monomials of a finite set $F=\{f_1,\ldots,f_k\}$ of generators $f_i$ with scope $s_i$ of $I$. Set $A=\bigcup_{i=1}^k \a_i+(\N^{s_i}\times 0^{n-s_i})$ and denote by $B=A^c =\N^n\sm A$ its complement. Write $K[[x]]^B$ for the space of power series whose expansions involve only monomials with exponent in $B$. An {\it echelon power series division} of a series $f\in K[[x]]$ by $F$ with respect to $<_\eps$ is a decomposition \med

\cl{$f = \sum_{i=1}^k a_if_i +b$}\med

with quotients $a_i\in K[[x_1,\ldots,x_{s_i}]]$, remainder $b\in K[[x]]^B$, and so that  \med

\cl{$\inin_\eps(f-b)=\min\, \{\inin_\eps(a_i)\cdot\inin_\eps(f_i),\, i=1,\ldots,k\}$,\hs 2cm (*)}\med

where the minimum refers to the ordering of the monomials of $K[[x]]$ induced by $<_\eps$. %The condition (*) will be essential in the proof of termination of the echelon standard basis algorithm, cf.~Thm.~2. [\hh As PC mentioned correctly, the original condition (*) without $\inin(b)$ would not have been sufficient. But one could even require in addition to (*) the condition that also $\inin(f')=\min\, \{\inin(a_i)\cdot\inin(f_i),\, i=1,\ldots,k\}$ where $f'$ denotes the part of the expansion of $f$ whose monomials belong to $\sum_{i=1}^k K[[x_1,\ldots,x_{s_i}]]\cdot\inin(f_i)$.]
This is the analog requirement as for polynomial division in the case where the divisors, say ideal generators, are not yet (or not necessarily) a Gr\"obner basis. Note that in general the decomposition is not unique. %We call $b$ {\it a} remainder of $f$ with respect to $f_1,\ldots,f_k$.

%{\it Example} : [silly, has nothing to do with echelons] Take $n=2$, $f_1=x^2 - xy$, $f_2= x^2y-y^3$, with $\inin(f_1)=x^2$, $\inin(f_2)=x^2y$, $s_1=s_2=2$, for which $f:=y^3-xy^2= 0\cdot f_1+0\cdot f_2+b=y\cdot f_1-f_2+0$ has two decompositions, with different remainders $b=y^3-xy^2$ and $b=0$.\med

\med
%-------------------------------------------------
%           THEOREM 1    DIVISION THEOREM
%-------------------------------------------------

{\bf Theorem 1.} (Division theorem for echelons) {\it Let $F=\{f_1,\ldots,f_k\}$ be a finite set of series $f_i$ in $K[[x_1,...,x_n]]$ with assigned scopes $0\leq s_i\leq n$. Choose a monomial order $<_\eps$ on $\N^n$. For every series $f$ there exists an echelon power series division (with respect to $<_\eps$)\med

\cl{$f=\sum_{i=1}^k\, a_if_i +b$}\med

of $f$ by $f_1,\ldots,f_k$ with respect to $<_\eps$.}\med

%-------------------------------------------------
%           REMARKS ON DIVISION THEOREM
%-------------------------------------------------

{\it Remarks.} (a) If $f_1,...,f_k$ form an echelon standard basis, the remainder $b$ of the division is unique (whereas the coefficients $a_i$ still need not be unique). Uniqueness of $b$ does not hold for arbitrary $f_1,...,f_k$. Prescribing support conditions on the coefficients $a_i$ as in the proof below by choosing a partition $A=\dot\cup A_i$ of the set $A$ and requiring $\supp(a_i\cdot\inin(f_i))\subset A_i$ for all $i$, both the coefficients $a_i$ and the remainder $b$ can be made unique (though they will depend on the chosen partition of $A$).\med

(b) If $I$ is the echelon generated by $f_1,...,f_k$ with scopes $s_1,...,s_k$, and if we assume that also $f$ belongs to $I$ and has assigned scope $s$, there is a natural way to assign to the remainder $b$, which then again belongs to $I$, a scope: namely, define it as the minimum of $s$ and the scopes $s_i$ for those $i=1,...,k$ for which $a_i\neq 0$. This value can either be maximized over all presentations $f=\sum_{i=1}^k\, a_if_i +b$ (finding the maximum value may not be constructive), or it can be made unique by choosing a partition $A=\dot\cup A_i$ and support conditions on the $a_i$ so that the presentation is unique.\med

(c) We can check by Thm.~1 effectively whether an element $f$ belongs to $I$ {\it up to degree} $d$, since then we only need a finite part of the echelon standard basis, namely those elements whose initial monomials are not larger than all degree $d$ monomials.\med

(d) With a little more work (using elementary Banach space techniques), the division statement of the theorem can be established for convergent power series, cf.~[HM, Thm.~5.1]. It does not hold in general for algebraic series.\med

(e) The division theorem can also be formulated for vectors of power series and finitary module echelons.\med

%-------------------------------------------------
%              PROOF OF DIVISION THEOREM
%-------------------------------------------------

{\it Proof.} We shall show that the $K$-linear map\med

\cl{$ u: \prod_{i=1}^k K[[x_1,\ldots,x_{s_i}]]\times K[[x]]^{A^c}\map K[[x]],$}\med

\cl{$ (a_1,\ldots,a_k,b)\map \sum_{i=1}^k a_if_i+b$}\med

is surjective. Along the way, we shall in addition show that every series $f\in K[[x]]$ has a preimage $(a,b)$ so that the condition $\inin(f-b)=\min\, \{\inin(a_i)\cdot\inin(f_i),\, i=1,\ldots,k\}$ holds.\med

Write $u=v+w$ where $v$ is the ``monomial approximation'' of $u$ given by the initial monomials $x^{\a_i}$ of $f_1,\ldots,f_k$, i.e., where $v$ is the linear map\med

\cl{$ v: \prod_{i=1}^k K[[x_1,\ldots,x_{s_i}]]\times K[[x]]^{A^c}\map K[[x]],$}\med

\cl{$ (a_1,\ldots,a_k,b)\map \sum_{i=1}^k a_i x^{\a_i}+b.$}\med

By definition of $K[[x]]^{A^c}$, the map $v$ is surjective. We shall choose a linear subspace $N$ of the first factor $\prod_{i=1}^k K[[x_1,\ldots,x_{s_i}]]$ so that the restriction $v_N$ of $v$ to $N\times K[[x]]^{A^c}$ becomes an isomorphism of $K$-vectorspaces. Using the inverse of $v_N$ we shall then show that also the restriction $u_N$ of $u$ is an isomorphism. From this the surjectivity of $u$ follows. %This allows us to define a scission of $v$ via an inverse of $v_N$ (a {\it scission} of $v$ is a $K$-linear map $\sigma: K[[x]]\map \prod_{i=1}^k K[[x_1,\ldots,x_{s_i}]]\times K[[x]]^{A^c}$ so that $v\circ\sigma\circ v =v$).
Our choice of $N$ will ensure in addition the requirement (*) in the decomposition $f=\sum_{i=1}^k a_if_i+b$. \med

To construct $N$, we proceed as in the classical case of Gr\"obner or standard bases by defining a suitable partition of the set of exponents $A=\bigcup_{i=1}^k \a_i+(\N^{s_i}\times 0^{n-s_i})$ [Gal, GP, HM]. %[We call such sets {\it combinatorial echelons} in $\N^n$, and the sets $\a_i+(\N^{s_i}\times 0^{n-s_i})$ slices or azulejos. These echelons are a priori assumed to be {\it finitary} because otherwise one would get arbitrary subsets of $\N^n$.]
There is no distinguished choice of the partition of $A$. Typically, one sets $A_1= \a_1+(\N^{s_1}\times 0^{n-s_1})$, and then defines $A_i=[\a_i+(\N^{s_i}\times 0^{n-s_i})]\sm \bigcup_{j<i}\, A_j$. Here its is advisable to order the elements $f_1,...,f_k$ by decreasing scope, $s_1\geq \ldots \geq s_k$, in order to exhaust $A$ first by larger regions.\med

%-------------------------------------------------
%    REDUNDANT TEXT ON SLICES AND PARTITIONS
%-------------------------------------------------

\ignore

[The following text until the continuation of the proof of the theorem is irrelevant for the rest of the note.] For practical issues, it is advisible here to start with a generator of maximal scope, i.e., to numerate the generators according to decreasing scope. Taking any numbering, there appears a subtle difference to the classical situation (where all scopes are equal to $n$), where the individual sets are ``connected''. For echelons, they are no longer connected, as is illustrated by the following \med

{\it Example} : Take $n=2$, $\a_1= (0,1)$ with scope $1$, $\a_2=(1,0)$ with scope $2$. Then $A_1= (0,1)+(\N\times 0)$ and $A\sm A_1 = (1,0)+(\N\times 0) \, \dot\cup\, (1,2)+\N^2$ is a disjoint union of the two ``slices'' $A_2=(1,0)+(\N\times 0)$ and $A_3=(1,2)+\N^2$. Accordingly, to get a partition, we have to enlarge our generator system to the elements $\a_1= (0,1)$ with scope $1$, $\a_2=(1,0)$ with scope $1$, and $\a_3=(1,2)$ with scope $2$, so that\med

\cl{$ A=A_1\,\dot\cup\, A_2\,\dot\cup\, A_3$.}\med

[In our paper about echelons from 2004, we discussed these partitions to some extent, we should look it up again to compare.]\med

Let $e=(e_1,\ldots,e_n)\in \{0,1\}^n$ be a vector with components $0$ and $1$. We denote by $\N^e$ the ideal of $\N^n$ of $n$-tuples with zero components whenever $e_i=0$.\med

%-------------------------------------------------

{\bf Lemma.} [?] {\it For any combinatorial echelon $A=\bigcup_i\, \a_i+(\N^{s_i}\times 0^{n-s_i})$ in $\N^n$ there exists a finite partition into slices of the form $S_j=\b_j+\N^e$, $e\in\{0,1\}^n$, which, moreover, is a refinement of the original decomposition of $A$.}\med

%-------------------------------------------------

{\it Proof.} [insert]\med

{\it Remarks.} (a) Further on we will not need that the slices $S_j$ are of the prescribed form. We could simply take any partition $A=\dot\cup\, A_i$ with $\a_i\in A_i$.  \med

(b) One cannot ensure that the slices are of the form $S_j=\b_j+(\N^{t_j}\times 0^{n-t_j})$, due to Seiler's example of the monomial ideal $I=(xy,yz,zx)$ which does not admit a direct sum decomposition into slices of the form $K[[x_1,\ldots,x_s]]\cdot x^\a$ [Se]. [End of irrelevant text.]\med

%-------------------------------------------------
%           CONTINUATION OF PROOF OF DIVISION THEOREM
%-------------------------------------------------

\recognize

%{\it Continuation of the proof of the theorem.}

So let us fix a partition $A=\dot\cup A_i$ of the set $A$. Denote by $N_i=K[[x]]^{A_i}\subset K[[x]]$ the subspace of series with exponents in $A_i$, and set $N=\prod_{i=1}^k N_i$. It is then clear that the restriction $v_N$ of $v$ to $N\times K[[x]]^{A^c}$ is an isomorphism of $K$-vectorspaces. Let $v_N\inv$ be its inverse. We prove that the restriction $u_N$ of $u$ to $N\times K[[x]]^{A^c}$ is also an isomorphism. For this it is sufficient to show that the composition $u\circ v_N\inv= (v+w)\circ v_N\inv = \Id_{K[[x]]}+w\circ v_N\inv$ is an isomorphism of $K[[x]]$. The ``formal inverse'' $\sum_{i=0}^\infty (-w\circ v_N\inv)^i$ defines a linear map from $K[[x]]$ to $K[[x]]$ since applying $w\circ v_N\inv$ to a power series $h$ increases its initial monomial with respect to the ordering of the monomials induced by $<_\eps$. It is therefore the inverse to $u\circ v_N\inv$. This shows that $u_N\circ v_N\inv$ and hence also $u_N$ are isomorphisms. It follows that the map $u$ is surjective as claimed.\med

It remains to show (*). We clearly have $\inin(f-b)\geq\min\, \{\inin(a_i)\cdot\inin(f_i),\, i=1,\ldots,k\}$. If strict inequality would hold, the equality  $f=\sum_{i=1}^k a_if_i+b$ would imply, because of $b\in K[[x]]^{A^c}$, that $\inin(a_i)\cdot\inin(f_i)=\inin(a_j)\cdot\inin(f_j)$ for some pair $i\neq j$. This is impossible since $A_i\cap A_j=\emptyset$. The theorem is proven.\qed

%-------------------------------------------------------
%            CONSTRUCTION OF STANDARD BASES
%-------------------------------------------------

We have already mentioned that echelon standard bases of echelons need no longer be finite. However, the ideas of Buchberger's algorithm apply as well to construct the elements one by one. This goes as follows.\med\med
\gb
%-------------------------------------------------
%            THEOREM 2 ECHELON STANDARD BASES
%-------------------------------------------------

{\bf Theorem 2.} (Echelon standard bases) {\it Let $F=\{f_1,\ldots,f_k\}$ be a set of power series with assigned scopes $0\leq s_i\leq n$ generating a finitary echelon $I$ in $K[[x_1,...,x_n]]$,\med

\cl{$I=\sum_{i=1}^k K[[x_1,...,x_{s_i}]]\cdot f_i$.}\med

Fix a monomial order $<_\eps$ on $\N^n$ satisfying condition $(\#)$. There exists an algorithm to enlarge $F$ iteratively so that, for every monomial $x^\a$ of the initial echelon $\inin(I)$ of $I$, one arrives after finitely many enlargements at a set $\wt F$ which contains an element $f$ of assigned scope $s$ in $I$ with $x^\a\in K[[x_1,\ldots,x_s]]\cdot \inin(f)$.}\med

%-------------------------------------------------
%           REMARKS
%-------------------------------------------------

{\it Remarks.} (a) Said differently, the algorithm produces in finitely many steps the elements of an echelon standard basis of $I$ up to any prescribed degree. An {\it enlargement} of $F$ is defined as a finite set $\wt F$ containing $F$ all whose elements belong again to $I$ and carry an assigned scope. In the proof, the algorithm for constructing these enlargements will be described explicitly.\med

(b) We do not pretend that the algorithm terminates in the sense that, in the construction of the echelon standard basis, after finitely many steps no more enlargements occur (and, in general, this will not happen). Moreover, even in the case where a finite echelon standard basis exists, the algorithm may produce infinitely many elements (most of which will be redundant). This is due to the fact that we do not apply division after each step. % (we have to renounce to division in order to keep control of the scopes).
\med

%-------------------------------------------------
%             PROOF
%-------------------------------------------------

{\it Proof.} We first explain the algorithm, and then show the required property. It is a variation of Buch\-berger's algorithm in the version of power series, with the additional requirement of respecting in each step the scopes of the involved elements. \med

Choose, for every pair $i\neq j$, the canonical minimal relation $(m_i,m_j)\in K[x]^2$ between the initial terms (i.e., initial monomials taken together with their coefficients) $e_i\cdot x^{\a_i}$ and $e_j\cdot x^{\a_j}$ of $f_i$ and $f_j$, where $e_i$ and $e_j$ denote the respective coefficients in $K$ and where $m_i$ and $m_j$ are terms with appropriate coefficients so that\med

\cl{$m_i\cdot e_i\cdot x^{\a_i}+m_j\cdot e_j\cdot x^{\a_j}=0$.}\med

It is clear that the relations are unique up to multiplication by constants in $K$. Set $S(f_i,f_j):=m_i f_i + m_j f_j$. These linear combinations of $f_i$ and $f_j$ satisfy $\inin(S(f_i,f_j))>_\eps \inin(m_if_i)= \inin(m_jf_j)$, i.e., there occurs a cancellation of (monomial multiples of) the initial monomials of $f_i$ and $f_j$.\med

In the algorithm, we will only consider linear combinations $g_{ij}: =S(f_i,f_j)$ for which both $m_i\in K[x_1,...,x_{s_i}]$ and $m_j\in K[x_1,...,x_{s_j}]$ respect the assigned scopes of $f_i$ and $f_j$. The other $S(f_i,f_j)$ will be discarded. Observe here that if $m_i$ or $m_j$ violate the scope condition then all monomial relations between $e_i\cdot x^{\a_i}$ and $e_j\cdot x^{\a_j}$ violate it.\med

We assign to the elements $g_{ij}$ thus obtained the scope $s_{ij}:=\min\{s_i,s_j\}$ and add them to the set $F$. This will be done with all pairs $(i,j)$ satisfying the scope condition. The resulting set $\wt F$ together with the assigned scopes of its elements is considered as the first enlargement of $F$. We then iterate the procedure with $\wt F$. This completes the description of the algorithm.\footnote
{${}^{(4)}$}{\smallTimes Observe here that in the subsequent enlargements one does not need to reconsider combinations $S(f_i,f_j)$ of elements $f_i,f_j$ which have been taken care of earlier.}\med

We now show that the algorithm fulfills the assertion of the theorem. Let  $x^\a$ be a monomial of $\inin(I)$. We have to prove that, after finitely many enlargements of $F$, there is an element $f\in F$ with assigned scope $s$ in $I$ so that $x^\a\in K[[x_1,\ldots,x_s]]\cdot \inin(f)$. \med

%To prove the claim, we will use condition (\#) on the monomial order. It is clear that the algorithm produces in each division step a remainder which is either $0$ or has an initial monomial not yet present among the initial monomials of the elements of the current set $F$. It follows that these new initial monomials leave eventually any given bounded subset of $\N^n$ as we iterate the algorithm. By condition (\#), they will become larger than $x^\a$ with respect to $<_\eps$ from a certain stage on.\med

%To prove this claim we assume, by the contrary, that there is a monomial $x^\a\in\inin(I)$ which never belongs to $K[[x_1,\ldots,x_s]]\cdot f\subset I$ while pursuing the algorithm, for any $f\in F$ of assigned scope $s$. This $x^\a$ would then be a monomial that is missed by the algorithm. \med

We may assume that we have already run the algorithm until arriving at a set $F=\{f_1,\ldots,f_k\}$ for which all subsequent {\it new} initial monomials appearing later in the algorithm are larger than $x^\a$. %(It is not clear whether we will need  this assumption later on.)
Indeed, $\inin(g_{ij})>_\eps \inin(m_if_i)= \inin(m_jf_j)$ is strictly larger than the maximum of $\inin(f_i)$ and $\inin(f_j)$. As we don't reconsider combinations $S(f_i,f_j)$ taken care of in earlier enlargements, it follows that the new initial monomials appearing after an enlargement are all larger than the minimum of the new initial monomials of the preceding enlargement. We conclude that the sequence of new initial monomials is unbounded. Hence, by hypothesis ($\#$) on the monomial order, the sequence must overtake $x^a$ eventually. \med

As $x^\a\in\inin(I)$ we may write $x^\a=\inin(\sum\, a_if_i)$ for some $a_i\in K[[x_1,\ldots,x_{s_i}]]$ and $f_i\in F$. % [It is irrelevant that the choice of the coefficients $a_i$ is not constructive since we just want to prove the claim.]
Set \med

\cl{$M:=\min\,\{\inin(a_i\cdot f_i),\, i=1,\ldots,k\}$.}\med

%-------------------------------------------------
\ignore

We may choose [this will not be mandatory for the rest of the proof] the coefficients $a_i$ so that the monomial\med

\cl{$M:=\min\,\{\inin(a_i\cdot f_i),\, i=1,\ldots,k\}$}\med

is maximized with respect to $<_\eps$ over all $a_i$ with $f=\sum\, a_if_i$: Indeed, the monomials $\min\,\{\inin(a_i)\cdot\inin(f_i),\, i=1,\ldots,k\}$ are bounded from above by $x^\a$; but $<_\eps$ allows no infinite bounded and strictly increasing sequences -- hence the maximum is attained.  \med

\recognize
%-------------------------------------------------

Clearly, $M\leq_\eps x^\a$. If $M=x^\a$ we are done: There is an $i$ so that $x^\a=M=\inin(a_i)\cdot\inin(f_i)$, hence $x^\a\in K[[x_1,\ldots,x_{s_i}]]\cdot \inin(f_i)$. \med

If $M<_\eps x^\a$, we will see that the algorithm enlarges $F$ to a set \med

\cl{$\wt F=\{f_1,...,f_k,f_{k+1},...,f_{\tilde k}\}$}\med

with assigned scopes $s_i$ for $f_i$, and we then construct a presentation $f=\sum_{i=1}^{\tilde k} \wt a_if_i$ of $f$ with coefficients $\wt a_i\in K[[x_1,...,x_{s_i}]]$ so that \med

\cl{$\wt M:=\min\,\{\inin(\wt a_i\cdot f_i),\, i=1,\ldots,\tilde k\}>_\eps\min\,\{\inin(a_i\cdot f_i),\, i=1,\ldots,k\}=M$.}\med

This procedure is then repeated. By hypothesis ($\#$) on the monomial order, the resulting strictly increasing sequence of monomials $M$, $\wt M$, ... must reach  $x^\a$ after finitely many iterations. That is what we want to prove.\med

%-------------------------------------------------

To do so, let $C$ be the set of indices $i$ with $\inin(a_i)\cdot\inin(f_i)=M$.  We necessarily have $\abs C\geq 2$, since, due to the inequality\med

\cl{$\inin(\sum_{i=1}^k a_if_i)>_\eps  \min\,\{\inin(a_i\cdot f_i),\, i=1,\ldots,k\}$,}\med

a cancellation of (monomial multiples of) initial monomials $\inin(f_i)=x^{\a_i}$ must occur in the sum $\sum_{i=1}^k a_if_i$. Let $c_i$ and $e_i$ in $K$ denote the coefficients of the monomials $\inin(a_i)$, respectively $\inin(f_i)$, of $a_i$, respectively $f_i$. It follows that the vector $r\in K[x]^k$ with entries $r_i=c_i\cdot \inin(a_i)$ if $ i\in C$ and $r_i=0$ otherwise forms a monomial relation in $\prod_{i=1}^k K[x_1,...,x_{s_i}]$ between the terms $e_i\cdot\inin(f_i)$, for $i=1,...,k$,\med

\cl{$\sum_{i=1}^k r_i\cdot e_i\cdot\inin(f_i)= \sum_{i\in C} c_i e_i\cdot \inin(a_i\cdot f_i)= (\sum_{i\in C} c_i e_i)\cdot M=0$.}\med

%[To easen the notation, we may also consider initial terms $\int(f_i)= e_i\cdot \inin(a_i)$ of series as the initial monomials multiplied with the respective coefficients, or define the initial monomials $\inin(f_i)$ from the very beginning including the coefficients.] The vector $r$ is thus a linear combination of the monomial generators [i.e., by def., their entries are terms] of the module echelon of relations between the initial terms $e_i\cdot \inin(f_i)$ of those $f_i$ for which $i\in C$. These monomial generators are considered as vectors in $K[x]^k$ by putting zeros at the irrelevant entries with index $i\not \in C$.\med

For each pair $j\neq \ell$ in $C$, denote by $m_{j\ell}\in K[x]^k$ the relation vector between the monomials $e_i\cdot \inin(f_i)$, $i=1,...,k$, whose only non-zero entries occur for indices $j$ and $\ell$ and are the terms $m_j$ and $m_\ell$ appearing in the minimal monomial relation $m_j\cdot e_j\cdot x^{\a_j}+m_\ell\cdot e_\ell\cdot x^{\a_\ell}=0$ defined earlier in the description of the algorithm,\med

\cl{$m_{j\ell}=(0,...,0,m_j,0,...,0,m_\ell,0,...,0)\in K[x]^k$.}\med

In order not to have to exclude the case $j=\ell$ we may set all $m_{jj}$ equal to $0$. We leave it as a (simple) combinatorial exercise to check that the vectors $m_{j\ell}'\in K[x]^{\abs C}$ obtained from $m_{j\ell}$ by taking only the components with index in $C$ form a generator system of the module echelon of relations between the terms $e_i\cdot \inin(f_i)$, for $i\in C$, respecting the scopes. \med

The entry of the vector $r$ at index $i$ belongs to $K[x_1,...,x_{s_i}]$, by the choice of $a_i$ in $K[[x_1,...,x_{s_i}]]$, and the same holds for $m_i\in K[x_1,...,x_{s_i}]$. Now notice that if we have a sum $h =\sum h_i$ with $\inin(h_i) = M$ for all $i$ and so that $\inin(h) > M$ then $h =\sum \lambda_{j\ell}\cdot S(h_j,h_\ell)$ for some $\lambda_{j\ell}\in K$. Furthermore, we have $\inin(S(h_j,h_\ell)) > M$ for all $j,\ell$. In view of this we may therefore write \med

\cl{$r=\sum_{j,\ell\in C} b_{j\ell}\cdot m_{j\ell}$,}\med

for some coefficients $b_{j\ell}$ which are monomials in $K[x_1,...,x_{s_{j\ell}}]$, with $s_{j\ell}=\min\{s_j,s_\ell\}$ as above. \med % [The coefficients $b_{j\ell}$ are not unique, and bad choices could produce $b_{j\ell}$ outside $K[x_1,...,x_{s_{j\ell}}]$. We may or should add an argument why our choice is possible.] \med

In the description of the algorithm we defined elements $g_{j\ell}=m_{j\ell}\cdot (f_1,...,f_k)= m_j f_j+m_\ell f_\ell$ with assigned scope $s_{j\ell}=\min\{s_j,s_\ell\}$ (the dot represents the scalar product in $K[[x]]^k$). We enlarge now $F$ to a set $\wt F$ by adding all $g_{j\ell}$, for $j,\ell\in C$. Denote by $f_{j\ell}=g_{j\ell}\in \wt F$ these new elements, and assign to them the scopes $s_{j\ell}:=\min\{s_j,s_\ell\}$. We get \med

\cl{$\sum_{i=1}^k r_if_i= \sum_{j,\ell\in C} b_{j\ell} (m_jf_j+m_\ell f_\ell)=\sum_{j,\ell\in C} b_{j\ell} g_{j\ell}=\sum_{j,\ell\in C} b_{j\ell} f_{j\ell}$.}\med

Decompose all $a_i$ into $a_i= c_i\cdot \inin(a_i)+a_i'$ for some $a_i'\in K[[x_1,...,x_{s_i}]]$ with $\inin(a_i')>_\eps\inin(a_i)$. Then \med

\hs 3.1cm $f=\sum_{i=1}^k a_if_i=\sum_{i\not\in C} a_if_i+\sum_{i\in C} (c_i\cdot \inin(a_i)+a_i')f_i$\med

\hs 3.4cm $=\sum_{i\not\in C} a_if_i+\sum_{i\in C} a_i'f_i+\sum_{i=1}^k r_if_i$\med

\hs 3.4cm $=\sum_{i\not\in C} a_if_i+\sum_{i\in C} a_i'f_i+\sum_{j,\ell\in C} b_{j\ell} f_{j\ell}$\med

\hs 3.4cm $=:\sum_{i=1}^k \wt a_i f_i +\sum_{j,\ell\in C}\wt a_{j\ell}f_{j\ell}$,\med

with respective coefficients $\wt a_i\in K[[x_1,\ldots,x_{s_i}]]$ and $\wt a_{j\ell}\in K[[x_1,\ldots,x_{s_{j\ell}}]]$, where $\wt a_i := a_i$ for $i\not \in C$ and $\wt a_i:=a_i'$ for $i\in C$, and where $\wt a_{j\ell}:=b_{j\ell}$. This is a new presentation of $f$ as a linear combination of elements of our enlarged set $\wt F$. The scopes are respected. The first summand in the last line satisfies by definition of $C$ and $a_i'$ the inequality\med

\cl{$\min\,\{\inin(\wt a_i)\cdot\inin(f_i),\, i=1,\ldots,k\}>_\eps M=\min\,\{\inin(a_i)\cdot\inin(f_i),\, i=1,\ldots,k\}$.}\med

As for the second summand, recall that $\inin(f_{j\ell})=\inin(g_{j\ell})>_\eps \inin(m_jf_j)=\inin(m_\ell f_\ell)$ and $\inin(m_j)=\inin(a_j)$ for all $j,\ell\in C$. This implies that also  \med

\cl{$\min\,\{\inin(\wt a_{j\ell})\cdot\inin(f_{j\ell}),\, j,\ell\in C\}>_\eps M$.}\med

We have found, after the enlargement of $F$ to $\wt F$, a presentation \med

\cl{$f=\sum_{i=1}^k \wt a_i f_i +\sum_{j,\ell\in C}\wt a_{j\ell}f_{j\ell}$}\med

of $f$ as a linear combination respecting the scopes of the elements of $\wt F$ and with larger value\med

\cl{$\wt M:=\min\,\{\inin(\wt a_i)\cdot\inin(f_i),\, \inin(\wt a_{j\ell})\cdot\inin(f_{j\ell});\, i=1,\ldots,k;\, j,\ell\in C\}$.}\med

Repeating the construction we produce by successive enlargements of $F$ a sequence of monomials $M<_\eps \wt M<_\eps \ldots$ which eventually attains $x^\a$. This is what had to be shown.\qed
\big\goodbreak

%-------------------------------------------------
%           CODE   CK JULY 31 2017
%-------------------------------------------------

\ignore

\rule[3pt]{\textwidth}{0.4pt}
\parbox[t]{0.15\textwidth}{INPUT:}
\parbox[t]{0.84\textwidth}{%
  $f_1,\dots,f_k\in K[[x_1,\dots,x_n]]$,\\
  $s_1,\dots,s_k\in\mathbb{Z}$ with $0\leq s_i\leq n$,\\
  monomial order $<_\eps$ on $\mathbb{N}^n$
  }\\[1ex]
\parbox[t]{0.15\textwidth}{OUTPUT:}
\parbox[t]{0.84\textwidth}{%
  standard basis $F$ of the echelon generated by $f_1,\dots,f_k$
  with scopes $s_1,\dots,s_k$, subject to $<_\eps$
  }\\[1ex]
\rule{\textwidth}{0.4pt}

\recognize

%-------------------------------------------------
%      CODE CK FEB 10 2018
%-------------------------------------------------

\def\small{}
\def\tt{}

{\bf Pseudo-code of algorithm of Theorem 2}\big

\hrule\medskip

INPUT: \hskip .45cm   $f_1,\dots,f_k\in K[[x_1,\dots,x_n]]$ with scopes $s_1,\dots,s_k\in\Z$, $0\leq s_i\leq n$,

\hskip 1.7cm  monomial order $<_\eps$ on $\N^n$,

\hskip 1.7cm  monomial $x^\a\in\inin(I)$, where $I$ is the echelon generated by $f_1,\dots,f_k$.

OUTPUT:  \hskip .1cm enlargement $F\supseteq\{f_1,\dots,f_k\}$ such that $F$ generates $I$ and there is an $f\in F$

\hskip 1.7cm  for which $x^\a\in K[[x_1,\ldots,x_s]]\cdot \inin(f)$ holds, where $s$ is the scope of~$f$.

\medskip\hrule\bigskip

{\small\tt 01: } $\ell := k$\par
{\small\tt 02: } $\ell_1 := 0$\par
{\small\tt 03: } $F := \{f_1,\dots,f_\ell\}$\par
{\small\tt 04: } {\bf while} $\ell_1 < \ell$ and $\neg\exists\;1\leq i\leq\ell: x^\a\in K[[x_1,\ldots,x_{s_i}]]\cdot \inin(f_i)$ {\bf do}\par
{\small\tt 05: } \qquad $P := \{\{i,j\}:\, 1\leq i<j\leq \ell \;\land\; j>\ell_1\}$\par
{\small\tt 06: } \qquad $\ell_1 := \ell$\par
{\small\tt 07: } \qquad {\bf for} $\{i,j\} \in P$ {\bf do}\par
{\small\tt 08: } \qquad\qquad compute minimal monomial relation $(m_i,m_j)$ of $f_i$ and $f_j$\par
{\small\tt 09: } \qquad\qquad {\bf if} $m_i\in K[x_1,\dots,x_{s_i}]$ and $m_j\in K[x_1,\dots,x_{s_j}]$ {\bf then}\par
{\small\tt 10: } \qquad\qquad\qquad $g_{ij} := m_if_i + m_jf_j$\par
{\small\tt 11: } \qquad\qquad\qquad {\bf if} $g_{ij} \not\in F$ or $\max\{s_u: f_u=g_{ij}\}<\min\,\{s_i,s_j\}$ {\bf then}\par
{\small\tt 12: } \qquad\qquad\qquad\qquad $\ell := \ell + 1$\par
{\small\tt 13: } \qquad\qquad\qquad\qquad $f_\ell := g_{ij}$\par
{\small\tt 14: } \qquad\qquad\qquad\qquad $s_\ell := \min\,\{s_i,s_j\}$\par
{\small\tt 15: } \qquad\qquad\qquad\qquad $F := F \cup \{f_\ell\}$\par
{\small\tt 16: } \qquad\qquad\qquad {\bf end if}\par
{\small\tt 17: } \qquad\qquad {\bf end if}\par
{\small\tt 18: } \qquad {\bf end do}\par
{\small\tt 19: } {\bf end do}\par
{\small\tt 20: } {\bf return} $F$\par
\bigskip

This algorithm, although it follows closely Buchberger's algorithm, does not
reduce the S-polynomials. If we wanted to include division with remainder into
the algorithm, its presentation would become much more complicated, which is
related to the determination of the scope of newly added elements. Clearly,
the scope of the new element should be the minimum of all scopes of elements
that were used in the division. But then, we have to record also intermediate
elements in the reduction with maximal possible scope. We illustrate the
problem with an example: assume that $f_1, f_2, f_3$ have the scopes $s_1 =
s_2 = 2$, and $s_3 = 1$. Then we assign to $g_{1,2} = m_1f_1 + m_2f_2$ the
scope~$2$, but after reducing it with~$f_3$, we have to assign scope~$1$. If
we only add the final result (with scope~$1$) to~$F$, we may hence miss an
element of the standard basis of~$I$. For actual computations, this conceptual
version of the algorithm may be very inefficient, and therefore, in the next
section, we will apply division with remainder to the S-polynomials, since the
above-mentioned problem does not occur there.

\bigskip
%-------------------------------------------------
%            BACK TO GABRIELOV'S EXAMPLE
%-------------------------------------------------

{\bf Analysis of Gabrielov's example}\med

We now return to the study of Gabrielov's example, where $f=1$ and $g=x\cdot(e^z-1)$ have assigned scope $x,y$, and where $h=yz-x$ has assigned scope $x,y,z$ (for clarity, we indicate instead of the value of the scope of the generators those variables which are allowed to appear in the series with which the generators are multiplied). As mentioned in the introduction, the explanation of the example does not require theorems 1 and 2, though these results help to put things in the right perspective.\med

Note first that the sums in our chosen echelons $I=\C\{x,y\}\cdot f+ \C\{x,y\}\cdot g + \C\{x,y,z\}\cdot h$ and $\wh I=\C[[x,y]]\cdot f+ \C[[x,y]]\cdot g + \C[[x,y,z]]\cdot h$ are direct: If we had a non-trivial linear relation \med

\hs 3cm $a(x,y)\cdot 1+b(x,y)\cdot x\cdot(e^z-1)+c(x,y,z)\cdot (yz-x)=0$,\med

setting $z={x\over y}$ would express $e^{x\over y}$ as a quotient of power series in $x$ and $y$, which is impossible. %\med
%
%-------------------------------------------------
%
Now order the monomials $x^i y^j z^k$ lexicographically by their exponents so that $z <_{lex} y <_{lex} x$. This order does not satisfy condition ($\#$) from above, since it allows bounded infinite strictly increasing sequences. This violation of ($\#$) does not alter the explanation of the example, and as the choice of $<_{lex}$ simplifies the presentation, we admit it here as an appropriate order.\med

The {\it initial monomials} of $f$, $g$ and $h$ with respect to $<_{lex}$ are $\inin(f) = 1$, $\inin(g) = xz$, and $\inin(h) = yz$. A lengthy check shows that the initial echelon $\inin(I)$ of $I$ equals\med

\hs 3cm $\ds \inin(I) = \C\{x,y\}\cdot 1+ \left(\sum_{k=1}^\infty \C\{x\}\cdot x^kz^k\right)\cap \C\{x,y,z\}+ \C\{x,y,z\}\cdot yz$.\med

The sum $\sum_{k=1}^\infty \C\{x\}\cdot x^kz^k$ is well defined as a subspace of $\C[[x,y,z]]$, since the degree of the summands tends to infinity. So we may take its intersection with $\C\{x,y,z\}$. Similarly, we have\med

\hs 3cm $\ds \inin(\wh I) = \C[[x,y]]\cdot 1+ \sum_{k=1}^\infty \C[[x]]\cdot x^kz^k+ \C[[x,y,z]]\cdot yz$.\med

Both subspaces are no longer {\it finitary} echelons since they require infinitely many ``generators''. We will not use these decompositions of $\inin(I)$ and $\inin(\wh I)$ in the sequel, but it is helpful to keep them in mind. \med

We now start the algorithm for the construction of the echelon standard basis of $I$. For our purposes, it will be convenient to take some shortcuts using Thm.~1 by dividing new elements by the preceding ones, so as to simplify the resulting series. Moreover, we will not show that the constructions produce eventually all initial monomials of $I$. In this sense, the analysis of the example relies on a slightly modified version of the algorithm of Thm.~2.\med

Recall that $f$, $g$ and $h$ have initial monomials $1$, $xz$ and $yz$. As $\C[[x,y]]\cdot 1 \cap \C[[x,y]]\cdot xz = 0$ and $\C[[x,y]]\cdot 1 \cap \C[[x,y,z]]\cdot yz = 0$, we can only take one linear combination, say, of $g$ and $h$, namely $S(g,h)=-y\cdot g +x\cdot h$. Division of this series by $f$, $g$ and $h$ with respect to the assigned scopes yields the element\med

\hs 3cm $g_2=-y\cdot g +x\cdot h + [x^2\cdot f+{1\over 2} x\cdot g+z\inv x \cdot (e^z-1-z)\cdot h]$\med

\hs 3.35cm ${}=x^2\cdot f-(y-{1\over 2} x)\cdot g+z\inv x \cdot (e^z-1)\cdot h$\med

\hs 3.35cm ${}={1\over 2}\cdot z\inv x^2\cdot (e^z\cdot (z-2) +z+2)$\med

\hs 3.35cm ${}= {1\over 2}\cdot x^2\cdot\sum_{k=2}^\infty \, {k-1\over (k+1)!}\cdot z^k$\med

\hs 3.35cm ${}={1\over 12}\cdot [x^2z^2 +{1\over 2}\cdot x^2z^3+{3\over 20}\cdot x^2z^4+{1\over 30}\cdot  x^2z^5+{1\over 168}\cdot x^2z^6+{1\over 1120}\cdot x^2z^7+\ldots]$.\med

It has assigned scope $x,y,$ and initial monomial \med

\hs 3cm $\inin(g_2)=x^2z^2.$ \med

This monomial does not belong to $\C[[x,y]]\cdot 1 +\C[[x,y]]\cdot xz + \C[[x,y,z]]\cdot yz$. So we have found a new initial monomial of $\inin(\wh I)$. We iterate the process of taking linear combinations and then reducing by division. The first few elements one obtains after $f$, $g$, $h$ and $g_2$ are \med\goodbreak

%\hs 3cm $g_2 ={1\over 12}\cdot [x^2z^2 +{1\over 2}\cdot x^2z^3+{3\over 20}\cdot x^2z^4+{1\over 30}\cdot  x^2z^5+{1\over 168}\cdot x^2z^6+\ldots]$,\med

\hs 3cm $g_3= {1\over 720}\cdot[x^3z^3+{1\over 2}\cdot x^3z^4 + {1\over 7}\cdot x^3z^5+{5\over 168}\cdot x^3z^6+{5\over 1008}\cdot x^3z^7+\ldots]$,\med

\hs 3cm $g_4= {1\over 100800}\cdot [x^4z^4 + {1\over 2} \cdot x^4z^5+ {5\over 36}\cdot x^4z^6+ {1\over 36}\cdot x^4z^7+\ldots]$,\med

\hs 3cm $g_5= {1\over 25401600}\cdot [x^5z^5 + {1\over 2} \cdot x^5z^6+{3\over 22}\cdot x^5z^7+{7\over 2461}\cdot x^5z^8+\ldots]$,\med

\hs 3cm $g_6= {1\over 10059033600}\cdot [x^6z^6 + {1\over 2} \cdot x^6z^7 + {7\over 52} \cdot x^6z^8 +{1\over 39} \cdot x^6z^9 + \ldots]$,\med

all with assigned scope $x,y$ and initial monomials of the form $x^kz^k$. We keep the coefficients in front of the brackets since they will play a crucial role later on.
%-------------------------------------------------
%
The general formula for these and the next elements $g_k$ appearing in the algorithm is\med

\hs 3cm $\ds g_k=x^k\cdot\sum_{i=k}^\infty q_{i,k}\cdot z^i$,\med

with coefficients $q_{i,k}$ given by\med

\hs 3cm $\ds q_{i,k}= {(i-1)! \over 4^{k-1}\cdot (i-k)!\cdot (i+k-1)! \cdot ({1\over 2})^{\ol {k-1}}}$,\med

where $({1\over 2})^{\ol {k-1}}$ denotes the rising factorial ${1\over 2}({1\over 2}+1)\cdots({1\over 2}+k-2)$.
%
%-------------------------------------------------
%
The expansion of $g_k$ results from the linear combination %(cf.~with Buchberger's {\it subtraction-polynomials}, nowadays known as $S$-polynomials)
\med

\hs 3cm $S(g_{k-1},h):=-y\cdot g_{k-1} +  q_{k-1,k-1}\cdot x^{k-1}z^{k-2}\cdot h$ \med

of $g_{k-1}$ and $h$ given by the relation $(-y, x^{k-1}z^{k-2})$ between their initial monomials $x^{k-1}z^{k-1}$ and $yz$, taking into account the factor $q_{k-1,k-1}$ in front of $g_{k-1}$. Then $g_k$ is obtained as the remainder of the division of $S(g_{k-1},h)$ by the series $f,g,h$ and $g_2,\ldots,g_{k-1}$ as described in Thm.~1. All $g_k$ are convergent series with assigned scope $x,y$. \med

%-------------------------------------------------

The vital observation here is that the coefficients $q_{k,k}$ of the initial monomials $x^kz^k$ of $g_k$ {\it tend very fast} to $0$: more precisely, the successive quotients $q_{k,k}/q_{k+1,k+1}$ are quadratic polynomials in $k$. This convergence is caused by the rapidly decreasing coefficients ${1\over k!}$ in $g=x\cdot(e^z-1)$ (most other transcendental series $g$ with this property would also produce counterexamples, whereas the phenomenon does not occur for algebraic power series.)\goodbreak \med

Rewrite now the series $g_k$ as linear combinations of the original generators $f,g,h$ of $I$,\med

\hs3cm $g_k=a_k\cdot f+b_k\cdot g+c_k\cdot h$,\med

with uniquely defined convergent series $a_k,b_k\in\C\{x,y\}$ and $c_k\in\C\{x,y,z\}$. They are given by the recursions \med

\hs 3cm $a_k = -y\cdot a_{k-1} + {1\over 4(2k-3)(2k-5)}\cdot x^2\cdot a_{k-2}$,\med

\hs 3cm
$b_k = -y\cdot b_{k-1} +  {1\over 4(2k-3)(2k-5)}\cdot x^2\cdot b_{k-2}$,\med

\hs 3cm
$c_k = -y\cdot c_{k-1} +  {1\over 4(2k-3)(2k-5)}\cdot x^2\cdot c_{k-2} -  z\inv\cdot(a_{k-1}\cdot f + b_{k-1}\cdot g + c_{k-1}\cdot h)$,\med

with $a_1 = 0$, $a_2 = x^2$, $b_1 = 1$, $b_2 = -y + {1\over 2}x$, $c_1 = 0$, $c_2 = -z\inv\cdot x\cdot (1 - e^z)$.
%
%-------------------------------------------------
%
The preceding formulas imply that $a_k, b_k$ are homogeneous polynomials in $x$ and $y$ of degree $k$, respectively $k-1$, while $c_k$ is a polynomial in $x,y,z, e^z, z\inv$ without poles. Note that in the expansions of $a_k$, $b_k$ and $c_k$ the monomials $x^2y^{k-2}$, $y^{k-1}$ and $xy^{k-2}$, respectively, appear with coefficients $\pm 1$. \med
%-------------------------------------------------

%[It would be nice to have here the closed formulas for the linear combinations expressing $g_k$ in terms of $f,g,h$. See the computations of CK from May 17 2017.]\med

\goodbreak
%-------------------------------------------------

The successive quotients\med

\hs 3cm $\ds {q_{k+1,k+1}\over q_{k,k}}= {1\over 4\cdot (2k+1)\cdot(2k-1)}$\med

of the coefficients of $g_k$ tend quadratically towards $0$. %, so that the  $q_{k,k}$ themselves tend very fast to $0$.
As\med

\hs 3cm $\ds {q_{i,k}\over q_{k,k}}= {k!\cdot (i-1)!\over (i-k)!\cdot (i+k-1)!}\leq 1$\med

for $i\geq k$, all coefficients $q_{i,k}$ of the series $g_k$ become comparatively small to $q_{k,k}$ while $k$ increases. This then implies that {\it infinite} linear combinations of the series $g_k$ with rapidly {\it increasing} coefficients may still produce {\it convergent} series. A typical example would be the convergent series\med

\hs 3cm $\ds e(x,z):= \sum_{k=2}^\infty  {1\over q_{k,k}}\cdot g_k(x,z)$.\med

Various other combinations of the series $g_k$ could be taken.
By construction, the series $e$ belongs to the intersection $\wh I\cap\C\{x,y,z\}$. We show that it does not belong to $I$. By uniqueness of the presentation, it suffices to write $e$ as a linear combination $e=a\cdot f+b\cdot g+c\cdot h$ of $f,g,h$ with divergent series $a,b,c$. Set $r_k= {1\over q_{k,k}}$ so that $e= \sum_{k=1}^\infty  r_k\cdot g_k$ and $a=\sum  r_k\cdot a_k$, $b=\sum r_k\cdot b_k$, $c=\sum  r_k\cdot c_k$ with $a_k$, $b_k$ and $c_k$ as defined above. As we noted earlier, the monomials $x^2y^{k-2}$, $y^{k-1}$ and $xy^{k-2}$ appear with coefficients $\pm 1$ in the expansions of $a_k$, $b_k$ and $c_k$, respectively. As the successive quotients $r_{k+1}/r_k$ tend quadratically with $k$ to infinity, it follows that the series $a,b,c$ diverge.\med

{\it Remark.} Here is an intuitive argument why Gabrielov's example works: The key is the replacement of the $y$-multiples $xyz^k$ of the monomials $xz^k$ of the series $g$ by $x^2z^{k-1}$, after multiplication of $g$ with $y$ and then applying the replacement of $y$ using $h=yz -x$. This construction corresponds to a shift by $(1,0,-1)$ of the exponents of $g$. The difference $g(x,y,z)- {x\over z}\cdot g(x,y,z)$ together with the iterates of this procedure creates a sequence of series whose initial monomials have coefficients tending rapidly to $0$. This, in turn, creates the explosion of the coefficients when expressing $e$ as a linear combination of $f$, $g$, $h$.\med

\med
%-------------------------------------------------
%        OUTLOOK
%-------------------------------------------------

{\bf Outlook}\med

Let us conclude with a remark on how to construct further counterexamples to linear nested Artin approximation for analytic functions. The three series $f = 1$, $g = x\cdot (e^z-1)$ and $h = yz - x$ proposed by Gabrielov have two key properties which make the example work: First, equating the last one to $0$, solving for $z$ (the exponent of $e^z$ in $g$), and substituting $z$ in $g$ produces $x\cdot(e^{x\over y} -1)$, with an essential singularity at $0$ with respect to $y$. This ensures that the sum in the echelon $I$ is direct, a fact which is needed to exclude non-trivial linear relations between $f$, $g$, and $h$. Secondly, the coefficients of $g$ tend sufficiently fast to $0$ (in the example, they are ${1\over k!}$). This allows to define a fourth {\it convergent} series $e$ as an infinite linear combination of the standard basis elements $g_k$ with rapidly {\it increasing} constant coefficients. As a consequence, the coefficient series $a$, $b$, $c$ in the presentation $e=a\cdot f+b\cdot g+c\cdot h$ must diverge. These are the only two properties used to produce the counterexample.\med

In view of this one could start off with an arbitrary $h=P(x,y)\cdot z-Q(x,y)$, with polynomials $P$ and $Q$ so that $Q/P$ has a pole at $0$, set again $f=1$ and take for $g$ any convergent power series in $z$ whose coefficients tend to $0$ at least as fast as ${1\over k!}$. Up to applying some simple algebraic modifications to $g$ and $h$ in order to avoid that the initial echelon $\inin(I)$ is finitely generated one will get again a counterexample.\med

Gabrielov's example and the preceding analysis of its underlying pattern is nicely  contrasted by the theorem of Eisenstein and Heine about the behaviour of the coefficients of algebraic power series [Eis, Hei, DP]: \med

{\it Let $h(x)\in \Q[[x]]$ be a univariate algebraic power series. There exists a positive integer $m$ so that $h(mx)\in\Z[[x]]$ has integer coefficients. In particular, the denominators of the coefficients of $h(x)$ have only finitely many prime factors and grow at most exponentially as $m^k$.}\med

This result may give a hint why there is no counterexample to Grothendieck's question for algebraic power series (a hint which is a fact by Popescu's theorem).

\big
%-------------------------------------------------
%\recognize
\med

{\bf References}\med

{\parindent 1cm \smallTimes \baselineskip 11pt

\litem{[ACH]} Alonso, M.E., Castro-Jim\'enez, F.J., Hauser, H.: Encoding algebraic power series. Found.~Comp.~Math. 2017, to appear.

%\litem{[Ab]} Abhyankar, S.: Two notes on formal power series. Proc. Amer. Math. Soc. 7 (1956), 903 - 905. %[Proves that $K[[x_1,x_2]]$ contains infinitely many formally independent series.]

%\litem{[AbM1]} Abhyankar, S., Moh, T.: On analytic independence. Trans. Amer. Math. Soc. 219 (1976), 77 - 87.

%\litem{[AbM2]} Abhyankar, S., Moh, T.: A reduction theorem for divergent power series. J. Reine Angew. Math. 241 (1970), 27 - 33.

%\litem{[AbP]} Abhyankar, S., van der Put, M.: Homomorphisms of analytic local rings. J. Reine Angew. Math. 242 (1970), 26 - 60.

%\litem{[ACH]} Alonso, M.E., Castro-Jim\'enez, F.J., Hauser, H.: Encoding algebraic power series. Manuscript 2015, 40 pp.

\litem {[Ar]} Artin, M.: Algebraic Spaces. Yale Univ. Press 1971.

%\litem {[Ar1]} Artin, M.: On the solutions of analytic equations. Invent. Math. 5 (1968), 277 - 291.

%\litem {[Ar2]} Artin, M.: Algebraic approximation of structures over complete local rings. Publ. Math. I.H.E.S. 36 (1969), 23 - 58.

%\litem{[BM]} [?] Bierstone, E., Milman, P.: Algebras of composite differentiable functions. In: Singularities, Part 1 (Arcata 1981), 127-136, Proc. Sympos. Pure Math. 40, Amer. Math. Soc. 1983.

\litem {[Bu]} Buchberger, B.: Ein algorithmisches Kriterium f\"ur die L\"osbarkeit eines algebraischen Gleichungssystems (An Algorithmic Criterion for the Solvability of Algebraic Systems of Equations). Aequationes math.~3 (1970), 374-383.

%Ein Algorithmus zum Auffinden der Basiselemente des Restklassenrings nach einem nulldimensionalen Polynomideal. Thesis, University of Innsbruck, Austria, 1966.

%\litem {[Cu]} Cutkosky, D.: Complex analytic morphisms and \'etoiles. Preprint 2015.

\litem{[DP]} Dwork, B., van der Poorten, A.: The Eisenstein constant. Duke Math.~J.~65 (1992), 23-43.

\litem{[Eis]} Eisenstein, G.: \"Uber eine allgemeine Eigenschaft der Reihen-Entwicklungen aller algebraischen Funktionen.
Bericht K\"onigl.~Preuss. Akad.~Wiss.~Berlin (1852), 441-443. Reproduced in Mathematische Gesammelte Werke, Band II, Chelsea Publishing 1975, pp.~765-767.

\litem{[Hei]} Heine, E.: Theorie der Kugelfunktionen, 2nd ed., Reimer Berlin, 1878.

\litem{[Gal]} Galligo, A.: A propos du Th\'eor\`eme de Pr\'eparation de
Weierstrass. Lecture Notes in Math. 409, 543-579. Springer 1973.

%\litem{[Ga]} Galligo, A.: Th\'eor\`eme de division et stabilit\'e en g\'eom\'etrie analytique locale. Ann.{\div} Inst.{\div} Fourier  39 (1979), 107-184.

\litem{[Gb1]} Gabrielov, A.: Formal relations among analytic functions. Funct. Anal. Appl. 5 (1971), 318-319.

\litem{[Gb2]} Gabrielov, A.: The formal relations between analytic functions.
Math. USSR. Izv. 7 (1973), 1056-1088.

\litem{[GP]} Greuel, G.-M., Pfister, G.: A Singular Introduction to Commutative Algebra, 2nd edition. Springer 2007.

%\litem{[Gra]} Grauert, H.: \"Uber die Deformation isolierter Singularit\"aten analytischer Mengen. Invent. Math. 15 (1972), 171 - 198.

\litem{[Gr]} Grothendieck, A.: S\'eminaire H. Cartan 1960/1961. Familles d'espaces complexes et fondements de la g\'eom\'etrie analytique. Fasc. 2, Exp. 13. Secr.~Math.~Paris 1962.

\litem{[HM]} Hauser, H., M\"uller, G.: A rank theorem for analytic maps between
power series spaces. Publ. Math.  I.H.E.S.  80
(1995), 95-115.

%\litem{[Hi]} Hironaka, H.: Idealistic exponents of singularity. In: Algebraic
%Geometry,  The Johns Hopkins Centennial Lectures. Johns Hopkins University Press 1977.

%\litem{[K+]} Kurke, H., Mostowski, T., Pfister, G., Popescu, D., Roczen, M.: Die Approximationseigenschaft lokaler Ringe. Lecture Notes Math. 634, Springer 1978.

\litem{[Os]} Osgood, W.: On functions of several complex variables, Trans. Amer. Math. Soc. 17 (1916), 1-8.

\litem {[Pa1]} Paw\l ucki, W.: On relations among analytic functions and geometry of sub-analytic sets. Bull.~Polish Acad.~Sci.~Math.~37 (1989), 117-125.

\litem {[Pa2]} Paw\l ucki, W.: On Gabrielov's regularity condition for analytic mappings. Duke Math.~J.~65 (1992), 299-311.

\litem{[Po1]} Popescu, D.: General N\'eron desingularization. Nagoya Math. J. 100 (1985), 97-126.

\litem{[Po2]} Popescu, D.: General N\'eron desingularization and approximation. Nagoya Math. J. 104 (1986), 85-115.

%\litem {[Se]} Seiler, W.: A combinatorial approach to involution and $\delta$-regularity I \& II.  Appl.{\div} Algebra Engrg.{\div} Comm.{\div} Comput.{\div} 20, no. 3-4 (2009), 207-259, 261-338.

\litem {[Sp]} Spivakovsky, M.: A new proof of D. Popescu's  theorem on smoothing of ring homomorphisms. J. AMS. 12 (1999), 381-444.

\litem {[Te]} Teissier, B.: R\'esultats r\'ecents sur l'approximation des morphismes en alg\`ebre commutative. %[d'apr\`es Artin, Popescu et Spivakovsky].
S\'em. Bourbaki 784 (1993/94).  Ast\'erisque 227 (1995), 259-282.

}

%-------------------------------------------------
\goodbreak\vs1cm

M.E.A.: Dept. de \'Algebra, \par
Universidad Complutense de Madrid, Spain\par
mariemi@mat.ucm.es\vs .5cm

F.J.C.-J.: Dept. de \'Algebra, \par
Universidad de Sevilla, Spain\par
castro@us.es\med

\vs -3.65cm \hs 7.5cm \vbox {H.H.: Faculty of Mathematics, \par
University of Vienna, Austria\par
herwig.hauser@univie.ac.at}\vs .6cm

\hs 7.5cm\vbox{C.K.: RICAM,\par
Austrian Academy of Sciences\par
christoph.koutschan@ricam.oeaw.ac.at}
%-------------------------------------------------

\vfill\eject\end